\documentclass{article}

\usepackage{arxiv}
\usepackage{amsmath}
\usepackage{amssymb}
\usepackage{amsthm}
\usepackage[utf8]{inputenc} 
\usepackage[T1]{fontenc}    
\usepackage{hyperref}       
\usepackage{url}            
\usepackage{booktabs}       
\usepackage{amsfonts}       
\usepackage{nicefrac}       
\usepackage{microtype}      
\usepackage{lipsum}		
\usepackage{graphicx}
\usepackage{natbib}
\usepackage{doi}
\newtheorem{theorem}{Theorem}

\newtheorem{corollary}[theorem]{Corollary}
\newtheorem{remark}[theorem]{Remark}
\usepackage{mathtools} 

\title{Non-Explosion Solutions for a Class of Stochastic
Physical Diffusion Oscillators}

\author{ \href{https://orcid.org/0000-0002-1972-379X}{\includegraphics[scale=0.06]{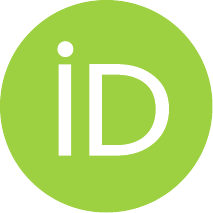}\hspace{1mm}Amrane  .~Houas}\\
	Laboratory of Applied Mathematics,\\  University Mohamed Khider, Biskra, Algeria.\\
	\texttt{houas.amrane@univ-biskra.dz} \\
	\And
	\href{https://orcid.org/0000-0002-6784-2216}{\includegraphics[scale=0.06]{orcid.pdf}\hspace{1mm} Fateh.~Merahi} \\
	Department of Mathematics\\
 University of Batna II, Algeria. \\
	\texttt{f.merahi@univ-batna2.dz} \\
 \And
	\href{https://orcid.org/0000-0002-8096-6280}{\includegraphics[scale=0.06]{orcid.pdf}\hspace{1mm} Mustafa.~Moumni} \\
	Laboratory of Photonic Physics and Nano-Materials,\\ University of Biskra, Algeria,\\
 LPRIM, Department of Physics, University of Batna I, Algeria.\\
	\texttt{m.moumni@univ-batna.dz} \\
}

\renewcommand{\shorttitle}{Non-Explosion Solutions for a Class of Stochastic
Physical Diffusion Oscillators}

\hypersetup{
pdftitle={Non-Explosion Solutions for a Class of Stochastic
Physical Diffusion Oscillators},
pdfauthor={Amrane.~Houas, Fateh.~Merahi, Mustafa. ~Moumni},
pdfkeywords={SDE, Non-explosion, Oscillator, Euler-Maruyama method}
}

\begin{document}
\maketitle

\begin{abstract}
	In this work, we are interested in problems that are related to the physical phenomena
of diffusion. We will focus on the theoretical aspect of the study, such as existence, uniqueness and
non-explosive solutions. We will weaken the conditions imposed on the coefficients of the stochastic
differential equations (SDE) that model some diffusion phenomena of mechanics. The work will be
based on a general non-explosion criterion and we will obtain sufficient conditions so that the solution
for a certain class of diffusions does not explode. We will construct Lyapunov functions that ensure
the non-explosion of the solutions. Two important oscillators, namely the Duffing and the Van Der Pol
oscillators, belong to this class. The Euler-Maruyama method is applied to these two oscillators to give
us a simulation solution for them.
\end{abstract}


\section{Introduction}
Differential equations for stochastic processes arise when studying many problems in physics and engineering (\cite{so13}, \cite{thygesen2023stochastic}, \cite{sarkka2019applied}, \cite{oksendal2022space}). They are typically of two basic types:
\begin{equation*}
\left\{
\overset{.}X_{t}=A(t)X_{t}+B(t) \hspace{2.5cm} \\
X_{t_{0}}=c \hspace{2.5cm} (I)
\right.
\end{equation*}
\begin{equation*}
\left\{
\overset{.}X_{t}=f(t,X_{t},\eta_{t}) \hspace{2.5cm} \\
X_{t_{0}}=c \hspace{2.5cm} (II)
\right.
\end{equation*}
where $ A(t) $, $ B(t) $ and $ \eta_{t} $ are random functions and the initial condition $ c $ is also a random variable in both cases.
If the random functions are sufficiently regular, we can consider the problems $ (I) $ and $ (II) $ as a family of problems corresponding to different realisations. In this case we can treat them with the classical methods of differential equations. The situation is different when the random functions are of the white noise type. White noise is perceived as a Gaussian stochastic process with zero mean and constant spectral density along the entire real axis \cite{ka12}. Such a process does not really exist because its covariance function is the Dirac function and therefore its variance is infinite. However, white noise is a very important mathematical idealisation for S.D.E. random influences that vary rapidly and whose values are virtually uncorrelated at different times \cite{ok03}.\\
The analysis of stochastic dynamical systems leads to differential equations of the following form
\begin{equation*}
\overset{.}{X}_{t}=f(t,X_{t})+G(t,X_{t})\xi_{t}
\end{equation*}
Here $ \xi_{t} $ is the White noise, $ G $ (respectively $ f $ ) is a function defined from $ [t_{0},t ]\times \mathbb{R}^{d} $ ( respectively $ [t_{0},t ]\times \mathbb{R}^{d} $ ) with values in $ Mat_{\mathbb{R}}(d,d) $ ( respectively $ \mathbb{R}^{d} $ ).\\
Some difficulties arise in studying these equations and they are due to the fact that white noise does not exist as a stochastic process in the classical sense. Other difficulties arise from the multitude of possible interpretations of this type of equation \cite{ik14}. \\
The challenges associated with this particular equation have been significantly mitigated due to the contributions of Ito, Bernstein, and Guikhman \cite{ik14}. Independently, the first two authors developed methods for addressing this type of equation (\cite{ito1951stochastic}, \cite{gikhman2007stochastic}). The solutions to these equations are Markovian processes with continuous trajectories and serve as probabilistic theoretical models for physical diffusion phenomena with random influences. A notable example of a diffusion process is the Wiener or Brownian motion process, as elaborated upon in the reference \cite{bu14}.\\
Non-linear stochastic oscillators have numerous important applications in various fields of science and engineering (\cite{ky15}, \cite{pik03}, \cite{br93}). This is because they correspond more closely to natural phenomena, as linearity is often just an approximation of reality that only serves to facilitate the theoretical study of real phenomena.
In this paper, our focus lies on the exploration of non-explosive solutions within the realm of a specific class of diffusion mechanics, which can be mathematically represented as::
\begin{equation*}
\overset{..}{x}_{t}+bf(x_{t},\overset{.}x_{t})+g(x_{t})=\sigma(x_{t},\overset{.}x_{t})\overset{.}{W}_{t}
\end{equation*}
Considering the solution as a vector in the phase space, the coefficients of drift and diffusion of the associated SDE do not satisfy the global conditions of existence and uniqueness of the solution. This implies that these solutions can escape to infinity in a finite time (known as an exploding SDE). To prevent this divergence, we establish sufficient conditions for $f$, $g$, and $\sigma$ based on non-explosion criteria in the general case \cite{kh11}.\\
Our focus is on the Duffing-Van der Pol oscillator (\cite{van1920theory}, \cite{van1927frequency}, \cite{duffing1918erzwungene}, \cite{guckenheimer2013nonlinear}), a special case of these oscillators. This oscillator plays a crucial role in various fields, including physics, chemistry, biology, electronics, and engineering  ( \cite{sz97}, \cite{ch11},\cite{bo03}, \cite{kudryashov2021generalized} )

\section{Motivation}
Modelling the dynamic behaviour of structures such as large suspension bridges and oil platforms, which are subjected to random oscillations caused by aerodynamic effects and earthquakes, leads to second-order equations of the following type:
\begin{equation}
\overset{..}{x}(t)+f(x(t),\overset{.}x(t))=\sigma(x(t),\overset{.}x(t))\overset{.}{W}(t)
\label{eq1}
\end{equation}
Where :
\begin{itemize}
\item[$ \bullet $] $ x(t) \in \mathbb{R}^{n}  $ with $ t \geq 0 $
\item[$ \bullet $] $ f $ is a function defined from $\mathbb{R}^{n}\times\mathbb{R}^{n}  $ to $ \mathbb{R}^{n} $
\item[$ \bullet $] $ \sigma $ is a function defined from $\mathbb{R}^{n}\times\mathbb{R}^{n}  $ to $ Mat_{\mathbb{R}}(n,m) $
\item[$ \bullet $] $ \overset{.}{W}(t) $  is a white noise of dimension $ m $
\end{itemize}
The first step consists of representing equation  (\ref{eq1}) as It\^{o} stochastic differential equation.\\
The classical method of approaching equations of type (\ref{eq1}) is to place the system in phase space by considering the process $ z(t)=(x(t),\overset{.}x(t)) $ and to represent (\ref{eq1}) as It\^{o} $ SDE $ on $ \mathbb{R}^{2n} $.\\
So if we set $y(t)= \overset{.}x(t)$, the equation (\ref{eq1}) is transformed to :
\begin{equation}
\left\{
\begin{aligned}
dx(t)&=y(t)dt \hspace{2.5cm} \\
dy(t)&=-f(x(t),y(t))+\sigma(x(t),y(t))dW(t) \hspace{2.5cm}
\end{aligned}
\right.
\label{eq2}
\end{equation}
The aim of this study is to investigate the existence and uniqueness of solutions to the equation ( \ref{eq2}) that do not explode in finite time.  \\
Due to the importance of Duffing-Van der Pol oscillators in the theory of stability and bifurcations for stochastic dynamical systems ( \cite{ba04}, \cite{ar93}, \cite{ar96}, \cite{ke99},  \cite{sc96},  \cite{li2023stochastic}, \cite{benedetti2022global}),  this study relies on two significant examples:\\
\textbf{Duffing oscillator} :
\begin{equation}
\overset{..}x(t)+2 \alpha\omega_{0}\overset{.}x(t)+\omega_{0}^{2}(x(t)+\lambda x(t)^{3})=\sigma\overset{.}W(t)
\label{eq3}
\end{equation}
\textbf{Van Der Pol oscillator}
\begin{equation}
\overset{..}x(t)+2 \xi\omega_{0}(x(t)^{2}-1)\overset{.}x(t)+\omega_{0}^{2}(x(t)+\gamma x(t)^{3})=\sigma\overset{.}W(t)
\label{eq4}
\end{equation}
\section{Non-Explosion Criteria}
Let's consider the SDE :
\begin{equation}
\left\{
\begin{aligned}
dx(t)&=b(x(t))dt+\sigma(x(t))dW(t), t>0 \hspace{2.5cm}\\
x(0)&=x_{0}     p.s \hspace{2.5cm}
\end{aligned}
\label{eq5}
\right.
\end{equation}
With 
\begin{itemize}
\item[$ \bullet $] $ x(t) \in \mathbb{R}^{n}  $ with $ t \geq 0 $
\item[$ \bullet $] $ b $ is a function defined from $\mathbb{R}^{n}  $ to $ \mathbb{R}^{n} $
\item[$ \bullet $] $ \sigma $ is a function defined from $\mathbb{R}^{n}  $ to $ Mat_{\mathbb{R}}(n,n) $
\end{itemize}

\begin{theorem}
(Non-Explosion Criteria) ( \cite{kh11} ) \\
 Under the hypothesis :
 \begin{itemize}
\item[\textbf{H1)}] $ \forall R>0 $, $ \exists K_{R}>0 $ : $ \vert b(x)-b(y) \vert + \vert \sigma(x)-\sigma(y) \vert \leqslant K_{R} \vert x-y \vert $,    $ \forall \vert x \vert \leqslant R $,  $ \forall \vert y \vert \leqslant R $
\item[\textbf{H2)}] $ x(0)  $ independent to $ W=\left\lbrace W(t), t \geqslant 0 \right\rbrace $ and $ E \vert x_{0} \vert^{2}<\infty $
\item[\textbf{H3)}] There exist a function $ x \mapsto V(x)  $ from $ \mathbb{R}^{n} $  to $ \mathbb{R}_{+} $ of $ C^{2} $ class, named the Lyapounov function and it verifies
\begin{equation*}
\left\{
\begin{aligned}
\forall x \in \mathbb{R}^{n} LV(x) \leqslant cV(x), c>0  \hspace{2.5cm} \\
V(x)\rightarrow +\infty if \vert x \vert \rightarrow +\infty \hspace{2.5cm}
\end{aligned}
\right.
\end{equation*}
With
\begin{equation*}
\left\{\begin{aligned}
L=\sum_{i=1}^{n}b_{i}(x)\dfrac{\partial}{\partial x_{i}}+\dfrac{1}{2}\sum_{i,j=1}^{n}a_{ij}(x)\dfrac{\partial^{2}}{\partial x_{i} \partial x_{j}}; \hspace{2.5cm} \\
a(x)=\sigma(x)\sigma^{T}(x) \hspace{2.5cm}
\end{aligned}
\right.
\end{equation*}
being the differential generator associated with the diffusion process solution of SDE
\end{itemize}
The SDE given by (\ref{eq5}) has a unique solution that does not explode in a finite time and it is a diffusion process.
\end{theorem}
\section{The main results}
\subsection{Diffusion of mechanics}
Consider the oscillator
\begin{equation}
\overset{..}x(t)+b(x(t),\overset{.}x(t))+g(x(t))=\sigma(x(t),\overset{.}x(t))\overset{.}W(t)
\label{eq6}
\end{equation}
Where
\begin{itemize}
\item[$ \bullet $] $ x(t) \in \mathbb{R}^{n} $, with $ t\geqslant 0 $
\item[$ \bullet$] $ b : \mathbb{R}^{2n} \rightarrow \mathbb{R}^{n}$
\item[$ \bullet$] $ g : \mathbb{R}^{n} \rightarrow \mathbb{R}^{n} $
\item[$\bullet$] $ \overset{.}W(t) $  is a white noise
\item[$\bullet $] $ \sigma : \mathbb{R}^{n}\times \mathbb{R}^{n} \rightarrow  Mat_{\mathbb{R}}(n,n)$
\end{itemize}
If we put $ y(t)=\overset{.}x(t) $, then equation \ref{eq6} will write as :
\begin{equation}
\left\{
\begin{aligned}
dx(t)&=y(t)dt \hspace{2.5cm}\\
dy(t)&= -[b(x(t),y(t))+g(x(t))]dt + \sigma(x(t),y(t))dW(t) \hspace{2.5cm}
\end{aligned}
\label{eq7}
\right.
\end{equation}
The differential generator associated with the equation (\ref{eq7}) is :
\begin{equation*}
L=<y,\nabla_{x}>-<b(x,y)+g(x),\nabla_{y} >+\dfrac{1}{2}Tr[\sigma(x,y)\sigma^{T}(x,y)D_{y}^{2}]
\end{equation*}
with
$ \nabla_{x}=(\dfrac{\partial}{\partial x_{1}},...,\dfrac{\partial}{\partial x_{n}}) $, $ \nabla_{y}=(\dfrac{\partial}{\partial y_{1}},...,\dfrac{\partial}{\partial y_{n}}) $ and $ D_{y}^{2}=((\dfrac{\partial}{\partial x_{i}\partial y_{i}})_{ij}) $, $ <.,.> $ is the scalar product
\begin{theorem} 
Let's consider the system
\begin{equation}
\left\{
\begin{aligned}
dx(t)&=y(t)dt \hspace{2.5cm}  \\
dy(t)&= -[b(x(t),y(t))+g(x(t))]dt + \sigma(x(t),y(t))dW(t) \hspace{2.5cm} \\
x(0)&=x_{0}, y(0)=y_{0} \hspace{2.5cm}
\end{aligned}
\label{eq8}
\right.
\end{equation}
Under the conditions :
\begin{itemize}
\item[\textbf{1)}] $ b $ is locally Lipschitzian,
\item[\textbf{2)}] The function $g$ derives from a potential i.e $ g(x)=\nabla_{x} G(x)  $ where $ G :\mathbb{R}^{n} \rightarrow \mathbb{R} $ of $ C^{2} $ class such that $G(x)\rightarrow +\infty  $ if $ \vert x \vert \rightarrow +\infty $
\item[\textbf{3)}] $ \exists K_{1} > 0  $, $ \exists c > 0 $ and 
\begin{equation*}
\dfrac{1}{2}Tr[\sigma(x,y)\sigma^{T}(x,y)] \leqslant <y,b(x,y)> + c[G(x)+\dfrac{1}{2} \vert y \vert^{2}]+K_{1}
\end{equation*}
\item[\textbf{4)}] The random variable $ (x_{0},y_{0}) $ satisfies the hypothesis \textbf{H2} of the theorem 1
\end{itemize}
The system \ref{eq8} admits a unique solution that does not explode in a finite time
\end{theorem}
\begin{proof} 
We will demonstrate that all the hypotheses of Theorem 1 are satisfied.\\
Firstly, it should be noted that $b$ is locally Lipschitzian. Secondly, $G$ is of $C^{2}$ class, which means that $g$ is of $C^{1}$ class. Additionally, considering condition 03, the coefficients of the system in \ref{eq8} are locally Lipschitzian ( \textbf{H1} is satisfied ).\\
Now consider the function $ V(x,y)=G(x,y)+\dfrac{1}{2}<y,y>+K $, where $ K $ is a positive constant that will be chosen further. Since $G(x)\rightarrow +\infty  $ if $ \vert x \vert \rightarrow +\infty $ , then $ \exists K_{2} > 0 : G(x)+K_{2} \geqslant 0  $ $ \forall x \in \mathbb{R}^{n} $ and so, we can choose the constant $ K $ sufficiently large to ensure that
\begin{equation*}
V(x,y) \geqslant 0 \forall (x,y) \in \mathbb{R}^{n} \times \mathbb{R}^{n}\\
V(x,y) \rightarrow + \infty if R= \vert x \vert^{2} + \vert y \vert^{2}  \rightarrow + \infty
\end{equation*}
Moreover, since $G$ is of $C^{2}$ class, $V$ is also of $C^{2}$ class.\\
Applying the differential generator associated with 8 to the function $V$ yields
 \begin{eqnarray*}
   LV(x,y) & = & <y, g(x) >-<y, b(x,y)+g(x) >+\dfrac{1}{2}Tr[\sigma(x,y)\sigma^{T}(x,y)] \\
  & = & -<y, b(x,y) > + \dfrac{1}{2}Tr[\sigma(x,y)\sigma^{T}(x,y)] \\
  & \leqslant &  c[G(x)+\dfrac{1}{2} \vert y \vert^{2}]+K_{1} \hspace{1.5cm} ( condition 03) ) \\
 & \leqslant & cV(x,y)-cK+K_{1}
 \end{eqnarray*}
We choose $ K \geqslant \dfrac{K_{1}}{c} $ such that $ LV(x,y) \leqslant cV(x,y) $, $ \forall (x,y)\in \mathbb{R}^{n} \times \mathbb{R}^{n} $\\
Therefore, the Hypothesis H3 of Theorem 1 is satisfied, which implies that the system (\ref{eq8}) has a unique solution that does not explode in a finite time.
 \end{proof}
 \begin{corollary}
 Under the conditions  1), 2) and 4) and if $ \sigma(x,y)= \sigma  $ is a constant and $ b $  verifies 
  \begin{center}
  $  \exists \alpha \geqslant 0   $ : $ <b(x,y),y> \geqslant -\alpha \vert y \vert^{2}  $ $ \forall (x,y) \in \mathbb{R}^{n}\times \mathbb{R}^{n} $,
\end{center}
the result of the theorem remains valid.
 \end{corollary}
 \begin{proof}
 We have\\
 $ [G(x)+\dfrac{1}{2} \vert y \vert^{2}] + <y, b(x,y) > + K \geqslant cG(x)+(\dfrac{c}{2}- \alpha )\vert y \vert^{2} + K $,  $ \forall (x,y) \in \mathbb{R}^{n}\times \mathbb{R}^{n} $.\\
Because $ Tr[\sigma(x,y)\sigma^{T}(x,y)] \geqslant 0 $, we have to choose the constant $ K$ sufficiently large and $ c > 2 \alpha $ such that $ cG(x)+(\dfrac{c}{2}- \alpha )\vert y \vert^{2} + K \geqslant \dfrac{1}{2}Tr[\sigma(x,y)\sigma^{T}(x,y)]$; hence the result follows.
\end{proof}
\subsection{Change of the variable}
The energy function of the system (\ref{eq8} ) serves as a Lyapunov function in Theorem 2. The function is defined as follows: $ V(x,y)=G(x,y)+\dfrac{1}{2}<y,y>+K $. However, it is necessary to assume that $G(x)\rightarrow +\infty  $ if $ \vert x \vert \rightarrow +\infty $; This is the first restriction. \\
A second restriction arises from the fact that if we apply the differential generator  $ L $ to this function, we have for $ y=0 $: $ LV(x,y)=\dfrac{1}{2}Tr[\sigma(x,y)\sigma^{T}(x,y)] $ and we cannot have $ LV(x,y) \rightarrow -\infty  $ if $ \vert x \vert^{2}+\vert y \vert^{2} \rightarrow + \infty $, which is a necessary condition to ensure the existence of a stationary solution.\\
To overcome these difficulties, we can apply a change of variables in the starting system and search for another Lyapunov function for the new equivalent system.\\
So, we assume that the system (\ref{eq6}) has the form:
\begin{equation}
\overset{..}x_{i}(t)+f_{i}(x_{i}(t))\overset{.}x_{i}(t)+\dfrac{\partial G}{\partial x_{i}}(x(t))=\overset{n}{\underset{j=1}{\sum }} \sigma_{ij}(x(t),\overset{.}x(t))\overset{.}W_{j} , i=\overline{1,n}
\label{eq9}
\end{equation}
where $ f_{i} $  $ (1 \leqslant i \leqslant n)$, is an application from  $ \mathbb{R} $ to $ \mathbb{R} $ and $ G $ is an application from  $ \mathbb{R}^{n} $ to $ \mathbb{R} $.\\
Let $ \overset{.}x(t)=y(t)=(y_{1}(t),...,y_{n}(t) $.\\
In phase space, the oscillator (\ref{eq9}) is written
\begin{equation}
\left\{
\begin{aligned}
  dx_{i}(t)=y_{i}(t) \hspace{1.5cm}\\
dy_{i}(t)=-[f_{i}(x_{i}(t))+\dfrac{\partial G}{\partial x_{i}}(x(t))]dt+ \overset{n}{\underset{j=1}{\sum }} \sigma_{ij}(x(t),y(t))\overset{.}W_{j} , i=\overline{1,n} \hspace{1.5cm}
\end{aligned}
\right.
\label{eq10}
\end{equation}
Now, consider the application
\begin{equation}
\left\{
\begin{aligned}
\Phi : \mathbb{R}^{n}\times\mathbb{R}^{n} \rightarrow \mathbb{R}^{n}\times\mathbb{R}^{n}\hspace{2.5cm} \\
(x,y)\mapsto (x_{1},...,x_{n},y_{1}+F(x_{1})),...,y_{n}+F(x_{n}))\hspace{2.5cm}
\end{aligned}
\label{eq11}
\right.
\end{equation}
where $ F_{i}(x_{i}) = \underset{0}{\overset{x_{i}}{\int }}f_{i}(s)ds $ \\
If we suppose that the functions $f_{i}$ are of class $ C^{1} $, then the application defined by (\ref{eq11}) is a $ C^{2} $ diffeomorphism from $ \mathbb{R}^{n} \times \mathbb{R}^{n} $ to  $ \mathbb{R}^{n} \times \mathbb{R}^{n} $.\\
Let's apply It\^{o}'s formula to the process
\begin{equation}
\left\{
\begin{aligned}
 \Phi_{i}(x,y)=x_{i} \hspace{3cm} \\
 \Phi_{i+n}(x,y)= y_{i}+F_{i}(x_{i}), i=\overline{1,n} \hspace{3cm}
\end{aligned}
\label{eq12}
\right.
\end{equation}
We obtain
\begin{equation}
\left\{
\begin{aligned}
d\Phi_{i}=(\Phi_{i+n}-F_{i}(\Phi_{i}))dt \hspace{2.5cm} \\
d\Phi_{i+n}=-\dfrac{\partial G}{\partial x_{i}}dt+\overset{n}{\underset{j=1}{\sum }}\sigma_{ij}(\varphi,\varphi -F(\varphi))dW_{j}(t) , i=\overline{1,n}\hspace{2.5cm}
\end{aligned}
\right.
\label{eq13}
\end{equation}
where $ \varphi = (\Phi_{1},...,\Phi_{n}) $, $ F(\varphi)=(F_{1}(\Phi_{1}),...,F_{n}(\Phi_{n})) $\\
Going back to the notation $ (x, y) $ in (\ref{eq13}), the system (\ref{eq10}) is therefore equivalent to the following system
\begin{equation}
\left\{
\begin{aligned}
dx_{i}=(y_{i}-F_{i}(x_{i}))dt \hspace{2.5cm} \\
dy_{i}=-\dfrac{\partial G}{\partial x_{i}}dt+\overset{n}{\underset{j=1}{\sum }}\sigma_{ij}(x,y-F(x))dW_{j}(t) , i=\overline{1,n} \hspace{2.5cm}
\end{aligned}
\right.
\label{eq14}
\end{equation}
where $ F(x)=(F_{1}(x_{1}),...,F_{n}(x_{n})) $\\
Let
\begin{equation}
H(x)= \overset{n}{\underset{i=1}{\sum }}\underset{0}{\overset{x_{i}}{\int }}F_{i}(s)ds
\label{eq15}
\end{equation}
then, the differential generator associated with the system (\ref{eq14}) is given by
\begin{equation}
L=<y-\nabla H(x), \nabla x>-<\nabla G(x), \nabla y>+\dfrac{1}{2} Tr[\sigma (x, y-F(x))\sigma^{T}(x, y-F(x))D_{y}^{2} ]
\label{eq16}
\end{equation}
Now we are interested in studying the non-explosive criteria for the solution of the system (\ref{eq14}) under weaker hypotheses than those of Theorem 2.
\subsubsection{Scalar case}
In this case, equation \ref{eq9} is written as
\begin{equation}
\overset{..}x(t)+f(x(t))\overset{.}x(t)+g(x(t))=\sigma(x(t),y(t))\overset{..}W(t)
\label{eq17}
\end{equation}
According to (\ref{eq14}), it is equivalent to the system
\begin{equation}
\left\{
\begin{aligned}
dx(t) = (y(t)-F(x(t)))dt \hspace{2.5cm} \\
dy(t) = -g(x(t))dt+\sigma(x(t),y(t)-F(x(t)))dW(t) \hspace{2.5cm}
\end{aligned}
\right.
\label{eq18}
\end{equation}
where $ F(x) = \underset{0}{\overset{x}{\int }}f(s)ds $
\begin{theorem}
Suppose that
\begin{itemize}
\item[\textbf{01)}] $ f $ and $ g $ are $ C^{1} $ class and $ \sigma $ is locally Lipschitzian,
\item[\textbf{02)}] $ \exists c_{1} > 0 $, $ \exists c_{2} > 0 $; $ x[F(x)+g(x)] \geqslant c_{1}x^{2} $, $ \forall \vert x \vert \geqslant c_{2} $,
\item[\textbf{03)}] $ \exists K_{1} > 0  $ : $ \dfrac{1}{2}\sigma^{2}(x,y-F(x)) \geqslant \dfrac{1}{2}(y-F(x))^{2}+F(x)[\dfrac{1}{2}F(x)+g(x)]+K_{1} $, $ \forall (x,y)\in \mathbb{R}^{2} $,
\end{itemize}
then for each initial condition that verifies the hypothesis \textbf{H2} of Theorem 1, the system (\ref{eq18}) has a unique solution that does not explode in a finite time.
\end{theorem}
\begin{remark}

This theorem shows that the following assumption
\begin{center}

$ \exists c_{1} > 0 $, $ \exists c_{2} > 0 $; $ xg(x) \geqslant c_{1}x^{2} $, $ \forall \vert x \vert \geqslant c_{2} $
\end{center}
previously thought to be necessary in precedent works (\cite{schenk1996deterministic}), is in fact not necessary.
\end{remark}
\begin{proof}
According to hypothesis \textbf{01}, the coefficients of system (\ref{eq18}) satisfy the condition of local Lipschitz.\\
Now, let us consider the function
\begin{equation*}
V(x,y)=\int_{0}^{x}[F(s)+g(s)]ds+\dfrac{1}{2}y^{2}+K
\end{equation*}
$ V $ is $ C^{2} $ class and the condition \textbf{02} implies that $ V(x,y )\geqslant 0 $ and $ V(x,y) \rightarrow +\infty $ if $ x^{2}+y^{2}\rightarrow +\infty $\\
Applying the differential generator associated with equation (\ref{eq18}) to $V$ yields:
\begin{eqnarray*}
LV(x,y)=(y-F(x))[F(x)+g(x)]-g(x)y+\dfrac{1}{2} \sigma^{2}( x,y-F(x)) \\
= yF(x)-F(x)[F(x)+g(x)]+\dfrac{1}{2} \sigma^{2}( x,y-F(x) \\
-\dfrac{1}{2}(y-F(x))^{2}-F(x)[\dfrac{1}{2}F(x)+g(x)]+\dfrac{1}{2}y^{2}+\dfrac{1}{2} \sigma^{2}( x,y-F(x))
\end{eqnarray*}
and using hypothesis \textbf{03} gives
\begin{equation*}
LV(x,y) \leqslant \dfrac{1}{2}y^{2}+K_{1} , \forall(x,y) \in \mathbb{R}^{2}
\end{equation*}
If we choose $ K \geqslant K_{1} $, then
\begin{equation*}
LV(x,y) \leqslant \dfrac{1}{2}y^{2}+K \leqslant V(x,y)
\end{equation*}
Therefore, this verifies the Hypothesis \textbf{H3} of Theorem 1, and as a result, the system (\ref{eq18}) has a unique solution that does not explode in a finite time.
\end{proof}

\subsubsection{Vectorial case}
\begin{theorem}
Assume that
\begin{itemize}
\item[\textbf{01}] $ f_{i} $, $ \dfrac{\partial G}{\partial x_{i}} $, are $ C^{1} $ class, for $ i=1,...,n $
\item[\textbf{02}] $ \underset{\vert x \vert \rightarrow +\infty} lim $  $  [H(x)+G(x)]=+\infty  $
\item[\textbf{03}] $ \exists K_{2} \geqslant 0 : Tr[\sigma(x,y-F(x))\sigma^{T}(x,y-F(x))] \leqslant <y,y>+<\nabla H(x),\nabla H(x)+2 \nabla G(x)> +K_{2}, \forall (x,y) \in \mathbb{R}^{n}\times \mathbb{R}^{n} $
\end{itemize}
then for each initial condition that verifies \textbf{H2} of Theorem 1, the system \ref{eq14} admits a unique solution that does not explode in a finite time
\end{theorem}
\begin{proof}
The condition \textbf{01} ensures that the hypothesis \textbf{H1} is verified.
Now consider the function
\begin{equation*}
V(x,y)=H(x)+\dfrac{1}{2}<y,y>+G(x)+K
\end{equation*}
where $ K>0 $ is a constant.\\
Because $ f_{i} $, $ \dfrac{\partial G}{\partial x_{i}} $ are $ C^{1} $ class for $ i=1,...,n $, then $ V $ is $ C^{2} $ class.\\
The condition \textbf{02)} implies that we have $ V(x,y)\geqslant 0 $, $ \forall (x,y) \in \mathbb{R}^{n}\times\mathbb{R}^{n}   $ and $ V(x,y) \rightarrow +\infty $ if $ R=\vert x\vert^{2}+  \vert y\vert^{2}\rightarrow +\infty $\\
Let us apply the differential generator associated with $ V$; we get:
\begin{eqnarray*}
LV(x,y)=<y-\nabla H, \nabla H +\nabla G> -<\nabla G,y>+\dfrac{1}{2}Tr[\sigma(x,y-F(x))\sigma_{T}(x,y-F(x))]\\
= V(x,y)-\dfrac{1}{2}<y,y>-\dfrac{1}{2}<\nabla H,\nabla H>- <\nabla H, \nabla G> + \dfrac{1}{2}Tr[\sigma(x,y-F(x))\sigma_{T}(x,y-F(x))]
\end{eqnarray*}
Using the condition \textbf{03)} and choosing $ K\geqslant K_{2} $ we get $ LV(x,y) \leqslant V(x,y) $, then the system (\ref{eq14}) admits a unique solution that does not explode at finite time.
\end{proof}
\section{Applications}
It will now be shown that both the Duffing and the Vander-Pol oscillators belong to the analyzed class of oscillators and that their solutions do not explode.
\subsection{Scalar case}
\textbf{Example 1 ( Duffing Oscillator)}\\
Let's consider the oscillator
\begin{equation}
\overset{..}x(t)+2 \alpha\omega_{0}\overset{.}x(t)+\omega_{0}^{2}(x(t)+\lambda x(t)^{3})=\sigma\overset{.}W(t)
\label{eq19}
\end{equation}
where $ \alpha $, $ \lambda $ and $\omega_{0}$ are constants in $\mathbb{R}_{+}^{*}  $ and $ \sigma $ is a constant. In this case, we have
\begin{equation*}
yb(x,y)=2\alpha\omega_{0}y^{2}, G(x)=\omega_{0}^{2}(\dfrac{x^{2}}{2}+\lambda\dfrac{x^{4}}{4})
\end{equation*}
then $ yb(x,y)\leqslant -y^{2} $ and $ G(x)\rightarrow + \infty $ if $ x\rightarrow + \infty $, therefore, corollary 1 applies in this case, which means that the solutions of the equation (\ref{eq19}) do not explode in a finite amount of time. \\
\textbf{Example 2 ( Vander-Pol Oscillator)}\\
We have the system:
\begin{equation}
\overset{..}x(t)+2 \xi\omega_{0}(x(t)^{2}-1)\overset{.}x(t)+\omega_{0}^{2}(x(t)+\gamma x(t)^{3})=\sigma\overset{.}W(t)
\label{eq20}
\end{equation}
where $ \xi $, $ \gamma $ and $\omega_{0}$ are constants in $\mathbb{R}_{+}^{*}  $ and $ \sigma $ is a constant. In this case, we have
\begin{equation*}
yb(x,y)=\xi\omega_{0}(x^{2}-1)y^{2}, G(x)=\omega_{0}^{2}(\dfrac{x^{2}}{2}+\gamma\dfrac{x^{4}}{4})
\end{equation*}
Here it is clear that corollary 1 is satisfied, meaning that solutions of (ref{eq20}) do not explode in finite time.\\
\textbf{Example 3 (Duffing-Van Der Pol Oscillator)}\\
Considering the oscillator
\begin{equation}
 \overset{..}x(t)+[\overset{2m}{\underset{j=1}{\sum }}\xi_{j}x(t)^{j}]\overset{.}x(t)+[\overset{2n}{\underset{j=1}{\sum }}a_{j}x(t)^{j+1}]=\sigma(x(t),\overset{.}x(t))\overset{.}W(t)
\label{eq21}
\end{equation}
with $ \xi_{2m} > 0$, $ a_{2n} < 0 $, and $\xi_{j} \in  \mathbb{R} $ for $ j=\overline{1,2m-1} $, $ a_{j } \in  \mathbb{R}$ for $ j=\overline{1,2n-1} $, and $ m > n \geqslant 1 $. The application $ \sigma $ satisfies the condition (03) of theorem 3.\\
In this case, we have:
\begin{center}
$ F(x)=\overset{2m}{\underset{j=1}{\sum }}\dfrac{\xi_{j}}{j+1}x^{j+1} $, $ g(x)=\overset{2n}{\underset{j=1}{\sum }}a_{j}x^{j+1} $ and
\end{center} 
\begin{center}
 $ \exists \alpha > 0, \exists  c_{1} \geqslant 0 : x[F(x)+G(x)] \geqslant \alpha x^{2}$, if $ \vert x \vert\geqslant c_{1} $.
 \end{center}
This verifies the conditions of Theorem 3, showing that solutions of (\ref{eq21}) do not explode in finite time.
\subsection{ Vectorial case}
\textbf{Example 1 ( Duffing Oscillator )}\\
We study the oscillator :
\begin{equation}
\overset{..}X(t)+B \overset{.}X(t)+g(X(t))=\sigma (X(t),\overset{.}X(t)) \overset{.}W(t)
\label{eq22}
\end{equation}
where $ X(t) \in \mathbb{R}^{n}$, $  \forall t \geqslant 0 $, $ B=(b_{ij}) $ is a positive matrix of rang $ n \times n $.\\
Suppose that
\begin{equation*}
g(x)=Ax+\overset{n}{\underset{i=1}{\sum }}K_{ii}x_{i}^{3}e_{i}, K_{ii} > 0,  \forall i=\overline{1,n}
\end{equation*}
where $ A=(a_{ij}) $ is a symmetric positive matrix of rang $ n \times n $ and $ \lbrace e_{i}\rbrace_{i=\overline{1,n}}  $ is the canonical base of $ \mathbb{R}^{n} $, then
\begin{equation*}
G(x)=\dfrac{1}{2}<x,Ax> +\dfrac{1}{4}\overset{n}{\underset{i=1}{\sum }}K_{ii}x_{i}^{4 } \geqslant \dfrac{1}{2}<x,Ax>\dfrac{\alpha}{2}\vert x \vert^{2}, ( \alpha > 0 ) ,
\end{equation*}
and so, $ G(x) \longrightarrow +\infty  $ if $ \vert x \vert \longrightarrow +\infty $\\
In addition we have
\begin{equation*}
<y,B(x,y)> = <y,By> \geqslant 0.
\end{equation*}
If the matrix $ \sigma(x,y) $ verifies :
\begin{equation*}
\exists K_{1} > 0, \exists \lambda \geqslant 0 : Tr[\sigma(x,y)\sigma^{T}(x,y)] \leqslant \lambda [G(x)+y^{2}]+K_{1}
\end{equation*}
then theorem 2 is applicable and the solutions of (\ref{eq22}) are non-explosive in a finite time.\\
\textbf{Example 2 }\\
Let be the oscillator :
\begin{equation}
\overset{..}x_{i}(t)+\xi_{i}(x_{i}^{2n_{1}}(t))\overset{.}x_{i}(t)+\dfrac{\partial G}{\partial x_{i}}(x(t))=\overset{n}{\underset{j=1}{\sum }}\sigma_{ij}(x_{i}(t),\overset{.}x_{i}(t))dW_{j}(t), i=\overline{1,n}
\label{eq23}
\end{equation}
with $ G(x)= - \overset{n}{\underset{j=1}{\sum }}(a_{i}x_{i}^{2n_{2}+2}+\nu x_{1} x_{i})$ where $ \xi_{i} > 0 $, $ a_{i} > 0 $ and $ \nu $ are real numbers and $ n_{1} > n_{2} > 0 $ are natural numbers.\\
In this case, we have
\begin{equation*}
f_{i}(x)=\xi_{i}x_{i}^{2n_{i}}, H(x)=\overset{n}{\underset{i=1}{\sum }}\dfrac{\xi_{i}}{(2n_{1}+1)(2n_{1}+2)}x_{i}^{2n_{i}+2}
\end{equation*}
The matrix $ \sigma(x,y) $ satisfies condition (03) of theorem 4 and, so, the solutions of system \ref{eq23} do not explode in a finite time.

\section{Simulation of Duffing and Vander-Pol oscillators with Euler-Maruyama method}
This section provides a simulation of the solution for Duffing and Vanderpol oscillators, expressed as a stochastic differential equation, using the Euler-Maruyama method.

The Euler-Maruyama method is a numerical technique used to simulate solutions to stochastic differential equations (SDEs). It is an extension of the classical Euler method for ordinary differential equations (ODEs) and provides an approximation of the trajectories of stochastic processes \cite{milstein2004stochastic}, particularly those involving random fluctuations like Brownian motion. This method is particularly useful when closed-form analytical solutions are difficult to obtain. The Euler-Maruyama method has various applications in different fields. The method is commonly used to study the dynamics of physical systems under random perturbations, simulate financial models with stochastic components, analyze biological processes influenced by noise, and more (\cite{higham2002strong}, \cite{lang2023euler}, \cite{wang2020euler}, \cite{jin2023central}, \cite{hiderah2023truncated}).
\subsection{\textbf{Simulation of Duffing oscillator}}
We have already shown the existence, uniqueness and non-explosivity of the Duffing oscillator solution in \ref{eq19} 
\begin{equation*}
\overset{..}x(t)+2 \alpha\omega_{0}\overset{.}x(t)+\omega_{0}^{2}(x(t)+\lambda x(t)^{3})=\sigma\overset{.}W(t)
\end{equation*}
where $ \alpha $, $ \lambda $, $\omega_{0}$, and $ \sigma $ are constants in $\mathbb{R}_{+}^{*}  $.\\
We are going to do a simulation of the solution for \eqref{eq19} by means of the Euler-Maruyama method.
\subsubsection{\textbf{Discretization steps with Euler-Maruyama method}}
Starting with the given SDE \eqref{eq19}, we perform the following steps
\begin{itemize}
    \item [1)] Select a small time step size $ \Delta t $ to discretize the continuous-time equation. This step size determines the granularity of the simulation.
    \item [2)] Set initial conditions for $x(0)$ and $v(o)$, where $v(t))=\overset{.}x(t)$. Also, set up arrays to store the values of $x$, $v$, and time $t$ at each step.
    \item [3)] Generate random numbers following a normal distribution, scaled by $\sqrt{\Delta t}$, to represent the increments of the Wiener process (Brownian motion) at each time step.
    \item [4)] Iterate over time steps from $t=0$ to $t=T$ (where $T$ is the total simulation time) with steps of $ \Delta t $. 
    \item [5)] For each step:
    \begin{itemize}
        \item [a)] Calculate the deterministic term $a$  and the stochastic term $b$ based on the SDE \eqref{eq19}
        \begin{align*}
            a= -2\alpha \omega v(t)-\omega^{2}x(t)-\omega^{2}\lambda x(t)^{3} \\
            b= \sigma \overset{.}W(t)
        \end{align*}
        \item [b)]Update $x$ and $v$ using the Euler-Maruyama update rule:
        \begin{align*}
            x(t+\Delta t)= x(t)+v(t)\Delta t
            \\
           v(t+\Delta t)= v(t)+a\Delta t+b
        \end{align*}
    \end{itemize}
  \item [6)]   After the simulation is complete, we plot the values of $x$ and $v$ over time $t$ to visualize the behavior of the Duffing oscillator.
\end{itemize}
Figure 1 shows the graph of the simulation solution of position and velocity for the following values of the Duffing Oscillator Parameters: $ \alpha=0.5 $, $ \lambda=3 $, $\omega_{0}=1$, and $ \sigma=2 $.
\begin{figure} [htb]
    \centering
    \includegraphics[scale=0.7]{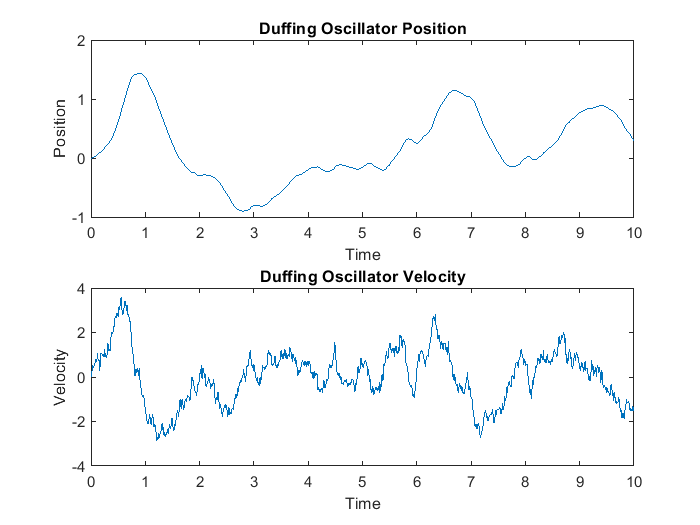}
    \caption{Duffing oscillator with Brownian movement: position and velocity}
    \label{fig:1}
\end{figure}

\subsection{\textbf{Simulation of Vander-Pol oscillator} }
In \eqref{eq20}, we have demonstrated the existence, uniqueness, and no-explosion of the solution of the Vander-Pol oscillator
\begin{equation*}
\overset{..}x(t)+2 \xi\omega_{0}(x(t)^{2}-1)\overset{.}x(t)+\omega_{0}^{2}(x(t)+\gamma x(t)^{3})=\sigma\overset{.}W(t)
\end{equation*}
where $ \xi $, $ \gamma $, $\omega_{0}$, and $ \sigma $ are constants in $\mathbb{R}_{+}^{*}  $
So, by means of the Euler-Maruyama method, we give a simulation of a solution for the Vander-Pol oscillator for some values of its parameters. 

\subsubsection{\textbf{Discretization steps with Euler-Maruyama method}  }
Starting with the given SDE \eqref{eq20}, we follow the same steps as for the Duffing case. Based on the equation \eqref{eq20} the deterministic and stochastic terms are given by:
\begin{align*}
    a= 2  \xi \omega x(t)^{2} v(t)- 2  \xi \omega v(t) + \omega^{2} x(t)
+ \omega^{2} \gamma x(t)^{3}  \\
b=  \sigma  \overset{.}{W(t)}
\end{align*}

Figure 2 shows the graph of the simulated solution for  $ \xi =0.1$, $ \gamma=0.25 $, $\omega=1$, and $ \sigma=0.1 $
\begin{figure}[htb!]
    \centering
    \includegraphics[scale=0.9]{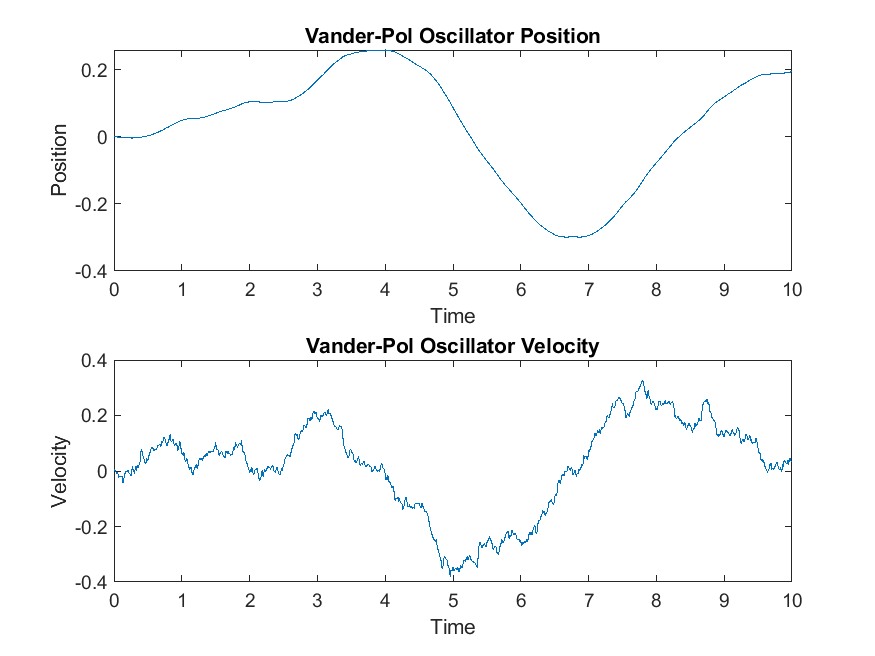}
    \caption{Vander-Pol oscillator with Brownian movement: position and velocity}
    \label{fig:2}
\end{figure}

\section{Conlusion}
This work focuses on the theoretical study of stochastic differential equations, particularly on the conditions for the existence, uniqueness, and non-explosion of solutions. The study is centred on problems related to the physical phenomena of diffusion.
The results show that it is possible to obtain sufficient non-explosive conditions. Sufficient conditions for constructing Lyapunov functions ensuring nonexplosive solutions are presented for a class of diffusion phenomena, including Duffing and Vander-Pol oscillators in scalar and vectorial cases. These reduced criteria has the potential to be applied to a wider range of physical cases. Finally, a simulation of Duffing and Vander-Pol oscillators using the Euler-Maruyama method is provided to demonstrate the non-explosive nature of the solutions.
\paragraph{\textbf{Declaration of Conflict of interest statement} }
All authors contributed to the study's conception and design. Material preparation, data collection, and analysis were performed by [Amrane Houas] and [Fateh Merahi], and [mustapha Moumni. The first draft of the manuscript to submission was written by [Amrane Houas] and the co-authors  [Fateh Merahi] and  [mustapha Moumni] commented on this version of the manuscript and approved it. The authors did not receive support from any organization for the submitted work.
\bibliographystyle{unsrtnat}
\bibliography{biblio}

\end{document}